\title{ An alternative proof that the Fibonacci group $F(2,9)$ is infinite.}
\author{ Derek F. Holt} 
\date{}
\documentstyle[12pt]{article}
\begin{document}
\maketitle
\begin{abstract}
This note contains a report of a proof by computer that the Fibonacci group
$F(2,9)$ is automatic. The automatic structure can be used to solve the
word problem in the group. Furthermore, it can be seen directly from the
word-acceptor that the group generators have infinite order, which of course
implies that the group itself is infinite.
\end{abstract}

The Fibonnacci groups $F(2,n)$, for integers $n \ge 2$, are defined by the
presentations
$$   \langle a_1, \ldots, a_n |
a_1a_2=a_3, a_2a_3=a_4, \ldots , a_{n-1}a_n=a_1, a_na_1=a_2 \rangle.
$$
They have been favourite test examples in combinatorial group theory for many
years. By 1974, it had been determined that they were finite for
$n = 2,3,4,5$ and $7$ and infinite for all other values of $n$, except possibly
$9$. The finiteness proofs were either by hand or by computer, using coset
enumeration, whereas the infiniteness proofs either constructed explicit
epimorphisms onto infinite groups or (for $n > 10$) used small
cancellation theory (see, for example, Sections 9 and 26 of [Johnson 1980]
for further details and references).
The remaining case $F(2,9)$
remained open until 1990, when it was finally proved infinite by Newman
in [Newman 1990]. In the meantime, it had been shown that this group had
a finite quotient of order $152 \cdot 5^{741}$ (see [Havas, Richardson
and Sterling 1979]), and Newman was
able to use the structure of this quotient, together with some
theoretical results, to show that it in fact had finite quotients of order
$152 \cdot 5^t$, for arbitrarily large values of $t$. This of course proved
that the group was infinite, but it provided no information about the
group itself, as distinct from its finite quotients. In particular, it
remained an open problem whether the generators $a_i$ had finite or infinite
order. (In fact, it still appears to be unknown whether there exists any
finitely presented infinite torsion group.)
In [Helling, Kim and Mennicke 1994], it is shown that
certain Fibonacci groups are fundamental groups of hyperbolic 3-manifolds,
which implies immediately that they are automatic (and in fact short-lex
automatic). However, these proofs do not apply to the groups $F(2,n)$ for odd
values of $n$.

The purpose of this note is to report on a successful computer proof that
$F(2,9)$ is automatic, and that its generators have infinite order.
For the general theory of automatic groups, see [Epstein et. al. 1992].
The algorithms employed in the computer programs that were used in the proof
are described in
[Epstein, Holt and Rees 1991] and, more recently, in [Holt 1994]. It was
much easier to prove the automaticity of $F(2,6), F(2,8)$ and
$F(2,10)$ computationally, since the associated automatic structures are much
smaller. The author had been attempting the calculation on $F(2,9)$ for
several years, but it was only recently that a computer with enough memory
became available.
Of course, with an enormous machine calculation of this nature, one
inevitably asks how far the result can be trusted, and whether it is likely
that a small logical or other error in the code (which cannot realistically be
ruled out as impossible) could have resulted in an incorrect final result.
I am of the opinion that this is extremely unlikely in this
particular case, for the following reasons. The correct automatic
structure is the result of a sequence of constructions of
approximations to this structure, and the eventual correct structure is
considerably smaller than the incorrect approximations. Furthermore, there is
a final stage to the computation, in which the structure is checked for
correctness by another program; this so-called ``axiom-checking'' program
constructs a series of much larger structures, in pairs, and the components of
each pair have to be identical for the verification process to succeed.
In addition, all of the calculations have been carried out successfully
by two radically different versions of the complete package,
and yielded the same results.

Roughly speaking, a group, with a given finite monoid-generating set $A$, is
automatic, if there exist two finite state automata $W$ and $M$, called the
{\em word-acceptor} and {\em multiplier}, respectively, with the following
properties.
The word-acceptor $W$ has input language $A$, and must accept a unique word in
$A^*$ for each group element. The multiplier $M$ reads pairs of words
$(w_1,w_2)$, with $w_1, w_2 \in A^*$, and accepts such a pair if and only if
$w_1$ and $w_2$ are both accepted by $W$ and $w_1^{-1}w_2$ is equal in $G$ to
one of the generators $a_i$. (See the references above for a more
detailed/accurate definition, and other equivalent definitions.) It turns out
that automaticity is a property of the group, and is independent of the
generating set, although the automata $W$ and $M$ will of course depend on
$A$. The Warwick programs that were used in the calculation for $F(2,9)$ are
only capable of calculating short-lex automtic structures, which are those in
which $W$ accepts the lexicographically least amongst the shortest words that
map onto a particular group element. (This assumes that an ordering has been
specified for $A$.) The automatic structure can be used to solve the
word-problem in $G$ efficiently, by reducing words in $A^*$ to their
representatives in the language of $W$.

For the calculation in $F(2,9)$, we used the ordered monoid generating set
$\{ a_1,a_1^{-1},a_2,a_2^{-1}, \ldots, a_9,a_9^{-1} \}$. We have not
experimented much with other generating sets or orderings, but the
ordering does not seem to have much influence on the difficulty of the
computation, whereas other generating sets (such as one with only two group
generators) seem to make things much more difficult. Our general experience
in this area suggests that the ``natural'' generating set is the best to use,
whenever this makes sense. Of course, the large number of generators does
increase the space requirements in some places. The final correct automata
$W$ and $M$ have 3251 and 25741 states, respectively. However, some of the
intermediate automata constructed were considerably larger than this;
further technical details follow below. 
By looking at the transitions of $W$, we observe immediately,
that when reading the word $a_1^n$ for $n>0$, $W$ starts in state 1, moves to
state 2, then to state 20, then to state 172 and then to state 686, where it
remains. Since all positive states of $W$ are accepting, this shows that $a^n$
is an accepted word for all $n$, which means that $a^n$ cannot be equal to the
identity for any $n>0$, and so $a$ has infinite order.

We conclude with some technical details of the computation.
This was done on a SPARCstation 20 with 256 Megabytes of core memory.
As mentioned above, it was actually successfully completed twice, the first
time using older versions of the software, which took several weeks of
continuous cpu-time. The second time, we were able to use a new version of
the code, which had been been completely rewritten in between. It then took a
total of about 12.5 hours of cpu-time, and used a maximum of just over 100
Megabytes of addressable memory, in addition to about 140 Megabytes of disk
space for temporary files. The files for the final correct automatic structure
have total size about 5 Megabytes.

The following description assumes that the reader has some familiarity with
the algorithm. The first step is to run the Knuth-Bendix process on the
presentation until the number of word-differences arising from reduction
equations appears to have become constant. This is probably the component of
the algorithm which has most scope for improvement, since the last few
word-differences seem to be very difficult to obtain.
We stopped with 539 word-differences (or 629 when closed under inversion),
at which point there were about 250000 reduction equations, and the process
was occupying about 100 Megabytes.  In fact, we did not have the full set of
word-differences at this point, which rendered the next few steps in the
calculation more difficult. The word-acceptor $W$ (calculated using all 629
word-differences) then had 8538 states. Using this and the word-differences
to calculate the multiplier $M$, resulted in $M$ having 1980342 states
initially, which minimized to 42808 states. It was this calculation which
required the large temporary filespace for storing the original unminimized
transition table for $M$ (this table does not need to be held in core memory,
but can be read in state-by-state during minimization).

The next step is a partial correctness test on $M$ (we test whether, for each
word $u$ accepted by $W$ and each generator $a_i$, there is a word $v$
accepted by $W$ such that $(u,v)$ is accepted by $M$, where $ua_i$ is equal to
$v$ in $G$). This test failed, and
increased the size of the (inverse-closed) word-difference set to 653.
We then recalculated
$W$ and $A$, which then had 8547 and 31021 states, respectively. The correctness
test failed again, and we then had 661 word-differences. This time, however,
$W$ had only 3251 states, and it turned out that $W$ was correct at this stage.
The number of states of $M$ was then 863871 before minimization and
25729 after minimization; the reduced sizes were due to the reduced size of
$W$. The correctness test failed twice more, but $W$ remained unchanged.
The number of word-differences increased to 671, and the number of states of
$M$ to 25741. At this stage, the partial correctness test succeeded, and we
could proceed to the full axiom-checking. This process took about 4.7 hours
cpu-time, and required about 105 Megabytes of core memory (which was the
largest amount of memory used at any stage).
The fact that all of the relations are short,
having length 2 or 3, rendered it more straightforward than usual, however.
From the correct automatic structure, it was then possible to construct
an automaton which accepts the minimal complete set of reduction rules,
and this in turn could be used to find the correct minimal set of
word-differences, of which there were 563. This is useful, since it can be used
to make the word-reduction process in the group more efficient.

\section{References}
[Newman 1990] M.F. Newman, ``Proving a group infinite'',
{\em Arch. Math.} {\bf 54} (1990), 209--211.\newline
\newline
[Epstein, Holt and Rees 1991] D.F. Holt, D.B.A. Epstein and S. Rees,
``The use of Knuth-Bendix
methods to solve the word problem in automatic groups'',
{\em J. Symbolic Computation} {\bf 12} (1991), 397--414.\newline
\newline
[Epstein et. al. 1992] D.B.A. Epstein, J.W. Cannon, D.F. Holt, S. Levy,
M.S. Paterson and W.P. Thurston,
``{\em Word Processing and Group Theory}'',
Jones and Bartlett, 1992.  \newline
\newline
[Havas, Richardson and Sterling 1979] G. Havas, J.S. Richardson and
Leon S. Sterling, ``The last of the Fibonacci groups'',
{\em Proc. Roy. Soc. Edinburgh} {\bf 83A} (1979), 199--203.\newline
\newline
[Helling, Kim and Mennicke 1994] H. Helling, A.C. Kim and J.L. Mennicke,
``On Fibonacci groups'', preprint.\newline
\newline
[Holt 1994] Derek F. Holt, ``The Warwick Automatic Groups Software'',
Submitted to Proceedings of DIMACS Conference on Computational Group Theory,
Rutgers, March 1994.\newline
\newline
[Johnson 1980] D.L. Johnson, ``{\em Topics in the Theory of Group
Presentations}'', London Math. Soc. Lecture Nore Series {\bf 42},
Cambridge University Press, 1980. \newline
\newline

\noindent Mathematics Institute,\\
University of Warwick,\\
Coventry CV4 7AL,\\
Great Britain.\\
e-mail: dfh@maths.warwick.ac.uk

\end{document}